\documentclass{article}
\usepackage{amsmath}
\usepackage{setspace}
\usepackage{amsthm}
\usepackage{amsfonts}

\newtheorem{Cor}{Corollary}
\newtheorem{Lem}{Lemma}
\newtheorem{Thm}{Theorem}
\usepackage{MnSymbol}
\usepackage{array}
\usepackage{cite}
\usepackage{hyperref}
\usepackage[T1]{fontenc}
\usepackage{ae}
\usepackage{aecompl}
\usepackage{bbm}
\usepackage{enumerate}
\mathchardef\mhyphen="2D
\author{Marc Fortier \footnote{mfortier0501@gmail.com}}
\title{Convergence Rates for Random Polarizations}

\begin{document}
\maketitle

\begin{abstract}
It is shown in \cite{AMF} that the expected $L^1$ distance between $f^*$ and $n$ random polarizations of an essentially bounded function $f$ with support in a ball of radius $L$ is bounded by $2dm(B_{2L})\|f\|_{\infty}n^{-1}$. This article expands on these results. We show that the expected $L^1$ distance is bounded by $c_n n^{-1}$ with $\limsup_{n\rightarrow \infty} c_n \leq 2^{d+1}\|\nabla f\|_1$ for every $f \in W^{1,1}(B_L) \cap L^{\infty}(B_L)$. Furthermore, we establish that the expected $L^1$ distance is $O(n^{-1/q})$ for $f \in L^p(B_L)$ with $1/p + 1/q = 1$. The rate $n^{-1}$ is shown to be best possible; specifically, $n$ times the measure of the symmetric difference between the random polarizations of a ball and its Schwarz symmetrization converges in distribution to a random variable with explicitly derived moments. We also prove that the expected symmetric difference between the random polarizations of a measurable set and its Schwarz symmetrization is slower than $n^{-r}$ for any $r > 2d$, and that if the rate is $n^{-1}$, the normalized symmetric difference converges in distribution. We introduce a new sequence of random polarizations where the transition probability depends on the state of the underlying Markov chain, yielding a convergence rate of $O\left( n^{-\left(2 - \frac{1}{d}\right)} (\log n)^{1 - \frac{1}{d}} \right)$ for $d>1$. Finally, we show that for every compact set $A \subset \mathbb{R}$ with finite perimeter, there exists a sequence of polarizations converging exponentially to its Schwarz symmetrization.
\end{abstract}

\section{Introduction}

It is well known that there exists sequences of polarizations (a rearrangement defined below) which can be applied iteratively to any initial function $f\in L^p$ $(1\leq p<\infty)$ to generate a sequence of functions (polarizations of $f$) which converge in $L^p$ to $f^*$ -- the symmetric decreasing rearrangement of $f$. The convergence is uniform when applied to continuous functions with compact support and, when applied to compact sets, the convergence also holds with respect to the Hausdorff distance (see \cite{AMF} for a detailed overview). The first result on rates of convergence of polarizations to the symmetric decreasing rearrangement appears in \cite[p.19]{AMF}: if $f\in L^1(\mathbb{R}^d)$ and bounded with support in $B_L$, then \begin{equation}\label{oldres}\mathbb{E}[||f^{\sigma_1\cdots \sigma_n}-f^*||_1] \leq 2dm(B_{2L})||f||_{\infty}n^{-1}.\end{equation}The purpose of this note is to expand on (\ref{oldres}). 

\section{Notation, rearrangements, results} 
\subsection{Notation}

In what follows, $m$ is the Lebesgue measure with sigma algebra $\mathcal{M}$; $B_{x,r}$ is the ball of radius $r$ centered at $x$; $B_r$ is the ball of radius $r$ centered at the origin; $\kappa_d$ is the volume of $B_1$ and $\omega_d$ is the surface area of $B_1$. The space of reflections that do not map the origin to the origin will be denoted by $\Omega$ and $\sigma_{x,y}$ will denote the unique reflection that maps $x$ to $y$.

\subsection{Rearrangements}
A rearrangement $T$ is a map $T:\mathcal{M}\rightarrow \mathcal{M}$ that is both \emph{monotone} ($A\subset B$ implies $T(A)\subset T(B)$) and \emph{measure preserving} ($m(T(A))=m(A)$ for all $A$). If \begin{equation}\mu_f(t) = m(\{f>t\})<\infty\end{equation} for all $t>0$, then $f$ is said to \emph{vanish at infinity} and we can define its rearrangement $Tf$ by using the ``layer cake principle": \begin{equation}\label{cocobee}Tf(x)=\int_{0}^{\infty}\mathbbm{1}_{T(\{f>t\})}(x)dt =\sup\{t: x\in T(\{f>t\})\}.\end{equation}

\subsubsection{Polarization}\label{go}

The \emph{polarization of $f$ with respect to $\sigma \in \Omega$} is defined as \begin{equation}\label{scotty1}f^{\sigma}(x)=\begin{cases}f(x)\vee f(\sigma(x)) & \text{if $x\in X_{+}^{\sigma}$} \\ f(x)\wedge f(\sigma(x)) & \text{if $x\in X_{-}^{\sigma}$} \\ f(x) & \text{if $x\in X^{\sigma}_{0}$}\end{cases}\end{equation}with $X^{\sigma}_{0}$ the hyperplane invariant under $\sigma$ which splits $\mathbb{R}^d$ into two disjoint half-spaces: $X^{\sigma}_{+}$, the half-space containing $0$, and $X^{\sigma}_{-}$, the half-space not containing $0$. If $A$ is an arbitrary set, then its polarization with respect to $\sigma$ is simply the polarization of $\mathbbm{1}_A$ and is denoted by $A^{\sigma}$: \begin{equation}\label{scotty}A^{\sigma}=\left(\sigma(A\cap X_{-}^{\sigma})\cap A^c \right)\cup \left(\sigma(A\cap X^{\sigma}_{+})\cap A\right) \cup \left(A\cap X_{+}^{\sigma}\right)\cup \left(A\cap X^{\sigma}_0\right). \end{equation}In other words, $A^{\sigma}$ is the same as $A$ except that the part of $A$ contained in $X_{-}^{\sigma}$ whose reflection does not lie in $A$ is replaced by its reflection in $X^{\sigma}_{+}$. As a result, polarization is measure preserving. It is clear from (\ref{scotty1}) that $f^{\sigma}\leq g^{\sigma}$ for all $\sigma \in \Omega$ whenever $f\leq g$ and thus polarization is monotone i.e., polarization is a rearrangement. One can check directly that $\{f^{\sigma}>t\}=\{f>t\}^{\sigma}$ for all $\sigma \in \Omega$ and, by (\ref{cocobee}), $f^{\sigma}$ is the rearrangement of $f$ with respect to the polarization rearrangement. 

\subsubsection{Schwarz Symmetrization}
For any $A\in \mathcal{M}$ there exists a unique open ball centered at the origin $A^*$ with the same measure as $A$ called the \emph{Schwarz symmetrization of $A$}. If $f(x)$ vanishes at infinity then its Schwarz rearrangement is denoted by $f^*(x)$. It is clear that $f^*(x)$ is radially decreasing: $f^*(x)\leq f^*(y)$ for $|x|\geq |y|$ and $f(x)=f(y)$ for $|x|=|y|$. In the literature, $f^*$ is also called the \emph{symmetric decreasing rearrangement of f}. If $f$ vanishes at infinity, we let \begin{equation}r_f(t)=\left(\frac{\mu_f(t)}{\kappa_d}\right)^{1/d}\end{equation} denote the radius of the open ball $\{f>t\}^*$. The distribution function $\mu_f(t)$ is always right continuous and thus so is $r_f(t)$. In particular, we have \begin{equation}\{f^*>t\}=\{f>t\}^*\end{equation} for all $t\geq 0$ and thus $f^*$ is right continuous.  

\subsection{Random Polarizations}
\subsubsection{Construction of probability measures}

Probability measures for $\Omega$ are easily constructed by mapping $\Omega$ to $\mathbb{R}^d$ via the invertible map $\varphi(\sigma) = \sigma(0)$: if $(\mathbb{R}^d,\mathcal{F},\mu)$ is a probability space, then $\mathbb{P}(\varphi^{-1}(A)) = \mu(A)$ yields a corresponding probability space for $\Omega$. A particularly good choice is $\mathcal{F}= \mathcal{M}$ and $d\mu = |x|^{-(d-1)} \mathbbm{1}_{B_{2L}} dm$: 
 \begin{equation}\label{eq:standard-prob}\mathbb{P}(\varphi^{-1}(A))=(2L\omega_d)^{-1}\int_{A \cap B_{2L}} \, |x|^{-(d-1)}dx\end{equation} for $A\in \mathcal{M}$. 
This probability measure has the following pleasant property (see \cite[p.18]{AMF}):
 \begin{equation}\label{chgprob}\mathbb{P}(\sigma(x)\in A)=(2L\omega_d)^{-1}\int_{A \cap B_{2L}}|x-y|^{-(d-1)}\, dy\end{equation}for every $x\in B_L$.

\subsubsection{Generating random polarizations}

We will be working with an i.i.d sequence of random polarizations $\sigma_n$: 

\begin{equation}\mathbb{P}(\bigcap_{i=1}^n\{\sigma_i \in \varphi^{-1}(B_i)\}) = \prod_{i=1}^n \mathbb{P}(\varphi^{-1}(B_i)).\end{equation} Given $f$ in $L^1$, we can iteratively apply the $n$ polarizations $\sigma_1,\ldots,\sigma_n$ to $f$ resulting in the random polarization $f^{\sigma_1\cdots \sigma_n}$. Since the $\sigma_n$ are independent, it is clear that $f^{\sigma_1\cdots \sigma_n}$ is a Markov chain. 

\subsection{Results}
\subsubsection{Rate of convergence estimates}
\begin{Thm}\label{bigtheorem}
Define $\rho_f(t) := \textnormal{per}(\{f > t\})$. The following holds:
\begin{enumerate}[(i)]

\item \label{mainprop:item:second}Suppose $f\in L^p(B_L)$ with $p \geq 1$ and let $x_n$ denote the unique fixed point of the function $dLn^{-1}(1+x^{\frac{1}{d}})^{d-1}$. If $p>1$ and $n$ is large enough that $x_n^{1/d} < r_f(0)$, then $\mathbb{E}[||f^{\sigma_1\cdots \sigma_n}-f^*||_1]$ is bounded by
\begin{equation}
\label{mainprop:item:second:eqn}
C_{d,p,L}||f||_p\;n^{-1}\left(\int_{L^{-1}x_n^{1/d}}^1 \frac{(1+r)^{(d-2)q}}{r^{(d-1)q/p}} dr\right)^{1/q} + 2||f^* \mathbbm{1}_{|x|\leq x_n^{1/d}}||_1
\end{equation}with $\frac{1}{p} + \frac{1}{q} = 1$ and $C_{d,p,L}=2(d-1)\omega_d^{1/q}L^{d/q + 1/p}$. If $p = 1$ and $n$ is large enough that $x_n^{1/d} < r_f(0)$, then $\mathbb{E}[||f^{\sigma_1\cdots \sigma_n}-f^*||_1]$ is bounded by \begin{equation}
\label{mainprop:item:second:eqn:2}
2(d-1)L^{d-1}||f||_1n^{-1}(1 + L^{-1}x_n^{1/d})^{d-2}x_n^{-(d-1)/d} + 2||f^* \mathbbm{1}_{|x|\leq x_n^{1/d}}||_1.
\end{equation}

\item \label{mainprop:item:third}If $m(\partial \{f > t\}) = 0$ for almost every $t$ and $f \in L^{\infty}(B_L)$, then \begin{equation}
\limsup_{n\rightarrow \infty}\; n \mathbb{E}[||f^{\sigma_1\cdots \sigma_n}-f^*||_1] \leq 2^{d+1}L \;||\rho_{f^*}||_1.\end{equation}

\item \label{mainprop:item:lbound} For any measurable set $A \subset B_L $, we must have $$\sum_{n=1}^{\infty}n\mathbb{E}[X_n]^{1/d} = \infty. $$ In particular, $n^r\mathbb{E}[X_n] \rightarrow \infty$ for $r>2d$.  

\item \label{mainprop:item:lbound-convdistr} If $\liminf_{n\rightarrow \infty}n\mathbb{E}[X_n] > 0 $ then $nX_n$ convereges in distribution. 

\end{enumerate}

\end{Thm}

(\ref{mainprop:item:third}) combined with the following two well-known properties of $W^{1,1}(B_L)$ yield the first corollary:

\begin{itemize}

\item \label{polya} If $f \in W^{1,1}(B_L)$, then, by the P\'{o}lya-Szeg\"{o} inequality, $f^* \in W^{1,1}(B_L)$ and $||\nabla f^*||_1 \leq ||\nabla f||_1.$
\item If $f \in W^{1,1}(B_L)$, then $||\rho_f||_1 = ||\nabla f||_1$. In particular, $m(\partial \{f > t\}) = 0$ for almost every $t > 0$. 

\end{itemize}

\begin{Cor}
\label{maincor:item:first}If $f\in W^{1,1}(B_L) \cap L^{\infty}(B_L)$, then \begin{equation}\limsup_{n\rightarrow \infty} \;n \mathbb{E}[||f^{\sigma_1\cdots \sigma_n}-f^*||_1] \leq 2^{d+1} L \;||\nabla f||_1.\end{equation}
\end{Cor}

\begin{Cor}
\label{maincor:item:second}If $f\in L^p(B_L)$ with $p>1$ and $\frac{1}{p} + \frac{1}{q} = 1$, then $\mathbb{E}[||f^{\sigma_1\cdots \sigma_n}-f^*||_1] = O(n^{-1/q})$.

\end{Cor}

\begin{proof}

If $f\in L^p(B_L)$ with $p>1$, then (\ref{mainprop:item:second:eqn}), H\"{o}lder's inequality and $x_n \sim dLn^{-1}$ gives \begin{align*}\mathbb{E}[||f^{\sigma_1\cdots \sigma_n}-f^*||_1] &= O(n^{-1}(\int_{x_n^{1/d}}^{L}r^{-(d-1)q/p} dr)^{1/q}) + O(||\mathbbm{1}_{|x|\leq x_n^{1/d}}||_q) \\
&= O(n^{-(1/q + 1/d)}) + O(n^{-1/q}) \\ &= O(n^{-1/q}).
\end{align*}

\end{proof}

\begin{Cor}
If $d=1$ then $nX_n$ converges in distribution.
\end{Cor}

\begin{proof}

It is shown in \cite{ABQD} that $X_n$ is always greater than or equal to $X'_n$ where $X'_n$ is the measure of the symmetric difference between the random polarizations of a non-centered ball, with the initial condition $X'_0 = X_0$, and its Schwarz symmetrization. Theorem (\ref{wcpb}) gives that $nX'_n$ converges in distribution, so $\liminf_{n\rightarrow \infty}nX_n > 0$. 

\end{proof}

\subsubsection{Random polarizations of balls}
We study the rate of convergence of the random polarization of balls (see \cite{AB} for a different approach). Let $A$ denote a ball of radius $r$ contained in $B_L$, $A_n = A^{\sigma_1 \cdots \sigma_n}$ and $X_n$ the distance from the origin of the centre of the ball $A_n$. 
\begin{Thm}\label{wcpb}
Let $u$ denote any unit vector. The following holds:
\begin{enumerate}[(i)]

\item \label{wcpb:first}The moments of $X_n$ can be computed exactly: \begin{equation}\label{exact2}\mathbb{E}[X_n^j] = \sum_{k=0}^{n}\binom{n}{k}(-1)^kX_0^{k+j}\prod_{i=0}^{k}c_{i+j-1}\end{equation} where $\label{exact}c_{\alpha}=(2L\omega_d)^{-1}\int_{|y|<1}(1-|y|^{\alpha})|u-y|^{-(d-1)}dy$ and $c_0 = 1$.

\item \label{wcpb:monotone}The moment $\mathbb{E}[X_n^k]$ is increasing in $X_0$ for all $k$ and all $n$. 

\item \label{wcpb:second}If $d=1$, then $\mathbb{E}[X_n] = 2L\int_0^{X_0/2L}(1-t)^n dt$ and $n X_n$ converges in distribution to an exponential distribution with scale parameter $2L$.

\item \label{wcpb:third} If $\alpha \geq 1$ then $\mathbb{E}[X_n^\alpha] \leq (X_0^{-1} + c_\alpha\alpha^{-1}n)^{-\alpha}$.
 
\item \label{wcpb:fourth}If $d>1$ then $\mathbb{E}[X_n] \geq (X_0^{-1} + \ell_dn)^{-1}$ where $$\ell_d = (2L\omega_d)^{-1}\int_{|y|<1}(1-|y|)|y|^{-1}|u-y|^{-(d-1)}dy.$$

\item \label{wcpb:sixth} $nX_n$ converges in distribution to a random variable $Y$ with moments $$\mathbb{E}[Y^k] = \frac{\alpha_d (2L)^k(k-1)!}{\prod_{i=1}^{k-1}{\overline{c}_i}}$$ where $$\alpha_d = \lim_{n\rightarrow \infty} n \sum_{k=0}^{n}\binom{n}{k}(-1)^k 2^{-(k+1)}\prod_{i=0}^{k}\overline{c}_i$$, $$\overline{c}_k = \omega_d^{-1}\int_{|y|<1}(1-|y|^k)|u-y|^{-(d-1)} \, dy $$, and $\overline{c}_0$ is set to 1.

\item \label{wcpb:fifth} $n \cdot m(A_n \triangle A^*)$ converges in distribution to $4\kappa_{d-1}r^{d} Y$. 

\end{enumerate}
\end{Thm}

\subsection{Proof of theorems}

We first consider a measurable set $A \subset B_L$ with $m(A) < m(B_L)$. Define the following random sequences:

\begin{itemize}

\item $A_n = A^{\sigma_1\cdots\sigma_n}.$
\item $X_n = m(A_n/A^*).$
\item $Y_n = (2L \omega_d)^{-1} m(A_n/A^{*})^{-2} \int_{A_n/A^*}\int_{A^*/A_n}|x-y|^{-(d-1)} dy dx.$
\item $r_n = (X_n/\kappa_d)^{1/d}.$ 

\end{itemize}

\begin{Lem}\label{mainlem}
The following holds: 
\begin{enumerate}[(i)]
\item \label{mainlem:item:first} $\mathbb{E}[X_n-X_{n+1} | X_{n}] = Y_nX_{n}^2.$
\item \label{mainlem:item:second} $\mathbb{E}[X_n^k-X^k_{n+1} | X_{n}] \geq Y_{n}X^k_{n} $ for all $k\geq 1.$
\item \label{mainlem:item:third} $\mathbb{E}[X_n^k] \leq \big(X_0^{-1} + k^{-1}\sum_{i=1}^n\mathbb{E}[Y_{i-1}^{-k}]^{-1/k} \big)^{-k}$.
\item \label{mainlem:item:fourth} $\mathbb{E}[X_n-X_{n+1} | X_{n}] \leq r_nL^{-1}X_{n}.$
\end{enumerate}
\end{Lem}

\begin{proof}
(\ref{mainlem:item:first})  Fubini's Theorem and (\ref{chgprob}) gives {\cite[p.19]{AMF}}:

\begin{align*}\mathbb{E}[X_{n}-X_{n+1} | X_{n}] &= \mathbb{E}[m(A_n/A^* \cap \sigma^{-1}(A^*/A_n)) ] \\ &=\int_{A_n/A^*} \mathbb{E}[\mathbbm{1}_{\sigma^{-1}(A^*/A_n)}(x)] dx \\ &= Y_{n}X_{n}^2. 
\end{align*}

\iffalse
(\ref{mainlem:item:second}) (\ref{mainlem:item:first}) gives \begin{align*}\mathbb{E}[X^2_{n} - X^2_{n+1}] &= 2\mathbb{E}[(X_{n}-X_{n+1})X_n] - \mathbb{E}[(X_n - X_{n+1})^2] \\ 
&\geq 2\mathbb{E}[X^3_nY_n] -  \mathbb{E}[m(A_n/A^* \cap \sigma^{-1}(A^*/A_n))^2] \\
&= 2\mathbb{E}[X^3_nY_n] - \mathbb{E}[\int_{A_n/A^*}\int_{A_n/A^*}\mathbbm{1}_{\sigma^{-1}(A^*/A_n)}(y)\mathbbm{1}_{\sigma^{-1}(A^*/A_n)}(z)dy dz] 
\end{align*}By (), we have  \begin{align*}\mathbb{E}[m(A_n/A^* \cap \sigma^{-1}(A^*/A_n))^2] &= \mathbb{E}[\int_{A_n/A^*}\int_{A_n/A^*}\mathbbm{1}_{\sigma^{-1}(A^*/A_n)}(y)\mathbbm{1}_{\sigma^{-1}(A^*/A_n)}(z)dy dz] \\&= \int_{A_n/A^*}\int_{A_n/A^*} \mathbb{P}(\{\sigma(y)\in A^*/A_n\} \cap \{\sigma(z)\in A^*/A_n\})dydz \\
&\leq X_n^3Y_n.
\end{align*}
\fi

(\ref{mainlem:item:second}) We proceed by induction. The case $k = 1$ is covered by (\ref{mainlem:item:first}). Suppose $k\geq 1$. We have \begin{align*}\mathbb{E}[X_n^{k+1} - X_{n+1}^{k+1}|X_n] &= \mathbb{E}[(X^k_{n}-X^k_{n+1})X_{n} - (X_{n+1}-X_{n})X^k_{n+1} |X_n] \\ 
&\geq \mathbb{E}[(X^k_{n}-X^k_{n+1})X_{n}|X_n] \\
&\geq Y_{n}X_n^{k+1}.
 \end{align*}

(\ref{mainlem:item:third}) By the mean value theorem and Jensen's inequality:

\begin{equation*}
\mathbb{E}[X^k_n]^{-1/k} - X_0^{-1} \geq k^{-1}\sum_{i=1}^n\frac{\mathbb{E}[X_{i-1}^{k+1}Y_{i-1}]}{\mathbb{E}[X^k_{i-1}]^{1 + \frac{1}{k}}}
 \geq k^{-1}\sum_{i=1}^n\mathbb{E}[Y_{i-1}^{-k}]^{-1/k}
\end{equation*} which is equivalent to
\begin{equation*}\mathbb{E}[X_n^k] \leq \big(X_0^{-1} + k^{-1}\sum_{i=1}^n\mathbb{E}[Y_{i-1}^{-k}]^{-1/k} \big)^{-k}
\end{equation*}

(\ref{mainlem:item:fourth}) Apply the Riesz rearrangement inequality:
\begin{align*}
\mathbb{E}[X_{n}-X_{n+1} | X_{n}] 
&\leq (2L\omega_d)^{-1}\int_{B_{r_n}}\int_{B_{r_n}}|x-y|^{-(d-1)} dy\; dx \\ 
&= r_n^{d+1}L^{-1}  (2\omega_d)^{-1}\int_{B_1}\int_{B_1}|x-y|^{-(d-1)}dy\; dx \\
&\leq X_nr_nL^{-1}.
\end{align*}

\end{proof}

We define the following random sequences associated with $f$:
\begin{itemize}
\item $A_{n,t} = \{f>t\}^{\sigma_1\cdots\sigma_n}$.
\item $X_{n,t} = m(A_{n,t}/A_{n,t}^*)$.
\item $Y_{n,t} = (2L \omega_d)^{-1} m(A_{n,t}/A_{n,t}^{*})^{-2} \int_{A_{n,t}/A^*}\int_{A^*/A_{n,t}}|x-y|^{-(d-1)} dy dx.$
\end{itemize}
To relate the rate of convergence for random polarizations of functions to that of measurable sets, we will make use of the following convenient formula:

\begin{equation}\label{layercakesymm}
||f^{\sigma_1 \cdots \sigma_n} - f^* ||_1 = \int_{0}^{||f||_{\infty}}m(A_{n,t} \triangle A_{n,t}^*) dt.
\end{equation}

\subsubsection{Proof of theorem \ref{bigtheorem}}

\begin{proof}

(\ref{layercakesymm}) and the previous lemma (\ref{mainlem}) gives
\begin{align}
\label{ubmain}\mathbb{E}[||f^{\sigma_1\cdots \sigma_n} - f^*||_1] &\leq 2\int_{0}^{b}(\sum_{i=1}^n\mathbb{E}[Y_{i-1,t}^{-1}]^{-1})^{-1}dt + 2\int_{b}^{\infty} \mu_f(t) dt \\
&\leq \label{ubmain2} 2\omega_d L n^{-1}\int_0^{b}(L + r_f(t))^{d-1}dt + 2\int_{b}^{\infty} \mu_f(t) dt
\end{align} for every $b \geq 0$.

(\ref{mainprop:item:second}) Now suppose that $f\in L^p(B_L)$ with $p\geq 1$ and $f^*(r)$ is absolutely continuous. Using the change of variable $t = f^*(r)$ and setting $b = f^*(x)$,  the right hand side of (\ref{ubmain2}) equals \begin{equation}\label{twoints}- 2L\omega_dn^{-1}\int_{x}^{r_f(0)}(L+ r)^{d-1} [f^*]'(r) dr - 2\kappa_d\int_{0}^{x} r^{d}[f^*]'(r)\;dr.\end{equation} Integration by parts shows that (\ref{twoints}) equals \begin{equation}\label{twoints-ip}\gamma(x;n) + 2L\omega_dn^{-1}(d-1)\int_{x}^{r_f(0)}(L+ r)^{d-2} f^*(r)dr + 2\omega_d\int_{0}^{x}f^*(r)r^{d-1}dr\end{equation} with $$\gamma(x;n) = 2\kappa_df^*(x)(dLn^{-1}(L+x)^{d-1} - x^d).$$ Recalling the sequence $x_n$ from statement (\ref{mainprop:item:second}), we have $\gamma(x^{1/d}_n;n) = 0.$ Assume that $p>1$. By H\"{o}lder's inequality ($\frac{1}{p} + \frac{1}{q} = 1$), the first integral in (\ref{twoints-ip}) is bounded by \begin{equation}
\label{ubint1} ||f||_p \;\omega_d^{-1/p}\left(\int_{x_n^{1/d
}}^L (L + r)^{(d-2)q}r^{-(d-1)q/p} dr\right)^{1/q}\end{equation} for $x=x_n^{1/d} < r_f(0)$. For $p = 1$, we have the upper bound \begin{equation}\label{ubint2}\int_{x_n^{1/d}}^{r_f(0)}(L+ r)^{d-2} f^*(r)dr \leq \omega_d^{-1}||f||_1(L + x_n^{1/d}
)^{d-2}x_n^{-(d-1)/d}.\end{equation}

(\ref{twoints-ip}) and (\ref{ubint1}) shows that the upper bound (\ref{mainprop:item:second:eqn}) is valid for all functions $f\in L^p(B_L)$ whose symmetric decreasing rearrangement $f^*$ is absolutely continuous. The set of all such functions is dense in $L^p(B_L)$: the P\'{o}lya-Szeg\"{o} inequality (\ref{polya}) implies that such a set must contain $W^{p,1}(B_L)$. Suppose that $f \in L^p(B_L)$, $f_k$ approaches $f \in L^p(B_L)$ with $f_k^*$ absolutely continuous, and $x_n^{1/d} < r_f(0)$. Convergence in $L^p$ implies that $$\lim_{k\rightarrow \infty} r_{f_k}(t) = r_f(t)$$ when $r_f$ is continuous at $t$. In particular,  $x_n^{1/d} < r_{f_k}(0)$ for sufficiently large $k$ and, for such $k$, the first inequality of (\ref{mainprop:item:second}) applies to $f_k$. Since both polarization and the symmetric decreasing rearrangement are contractive on $L^p$, \begin{align}
||(f_k^{\sigma_1\cdots \sigma_n} - f_k^*) - (f^{\sigma_1\cdots \sigma_n} - f^*)||_1 &\leq ||f^* - f_k^*||_1 + ||f-f_k||_1 \\
&\leq 2||f-f_k||_1
\end{align} for every random sequence $\sigma_n$. By the bounded convergence theorem: $$\lim_{k\rightarrow \infty}\mathbb{E}[||f_k^{\sigma_1\cdots \sigma_n} - f_k^*||_1] = \mathbb{E}[||f^{\sigma_1\cdots \sigma_n} - f^*||_1].$$ Since $||f_k^* \mathbbm{1}_{|x| \leq x_n^{1/d}}||_1$ tends to $||f^* \mathbbm{1}_{|x| \leq x_n^{1/d}}||_1$ then the first inequality in (\ref{mainprop:item:second}) holds for $f$. The same approximation procedure can be used to extend the second inequality in (\ref{mainprop:item:second}) to all of $L^1(B_L)$.

(\ref{mainprop:item:third}) We suppose that $||f||_{\infty} < \infty$ and $m(\partial \{f>t\}) = 0$ for almost every $t$. Recalling (\ref{ubmain}): 
\begin{equation}
\mathbb{E}[||f^{\sigma_1\cdots \sigma_n} - f^*||_{L_1}] \leq 2\int_{0}^{||f||_{\infty}}(\sum_{i=1}^n\mathbb{E}[Y_{i-1,t}^{-1}]^{-1})^{-1} \,dt.  
\end{equation} We have $$Y_{n,t} \geq (2\omega_d L)^{-1}[2r_f(t) + \underset{x \in A_{n,t}}{\sup} d_H(x,A_{n,t}^*)]^{-(d-1)}.$$

If $m(\partial \{f>t\}) = 0$ then $$\underset{x \in A_{n,t}}{\sup} d_H(x,A_{n,t}^*) \rightarrow 0$$ almost surely. Hence, by the dominated convergence theorem, 

$$\limsup_{n\rightarrow \infty} n \mathbb{E}[||f^{\sigma_1\cdots \sigma_n} - f^*||_1] \leq 2^{d+1}L \;||\rho_{f^*}||_{1}.$$

(\ref{mainprop:item:lbound}) We have

$$ 
\begin{aligned}
\mathbb{P}(X_{n+1} \leq X_n (1 - 1/n) ) &= \mathbb{E}[\mathbb{P}(X_{n+1} \leq X_n (1 - 1/n) | X_n)] \\
& \leq n L^{-1} \kappa_d^{-1/d} \mathbb{E}[X_n^{-1} X_n^{1+1/d}] \\
& \leq n L^{-1} \kappa_d^{-1/d}  \mathbb{E}[X_n]^{1/d}.
\end{aligned}
$$

By the Borel-Cantelli Lemma, if $$\sum_{n=1}^{\infty}n\mathbb{E}[X_n]^{1/d} < \infty $$ then $$X_{n+1} \geq (1-1/n)X_n $$ for $n \geq N$ (where the $N$ depends on the sequence), almost surely. But this would imply that $$\liminf_{n \rightarrow \infty}nX_n > 0$$ almost surely, and, by Fatou's Lemma, $$\liminf_{n \rightarrow \infty} n\mathbb{E}[X_n] > 0 $$, which is a contradiction.

(\ref{mainprop:item:lbound-convdistr}) We first show that $n^k \mathbb{E}[X_n^k]$ converges to a non-zero limit for all $k$. We have 

$$0 < \liminf_{n \rightarrow \infty}n^k \mathbb{E}[X_n]^k\leq \liminf_{n \rightarrow \infty}n^k \mathbb{E}[X_n^k] \leq \limsup_{n \rightarrow \infty}n^k \mathbb{E}[X_n^k] < \infty.$$ Hence it suffices to prove that 

$$\lim_{n \rightarrow \infty} \frac{\mathbb{E}[X_n^k - X^k_{n+1}]}{\mathbb{E}[X_n^k]} = 0.$$

We have 

$$
\begin{aligned}
\mathbb{E}[\Delta X_n^k] &\leq k \mathbb{E}[\Delta X_n X_n^{k-1}] \\
&\leq \mathbb{E}[(\Delta X_n)^2]^{1/2}\mathbb{E}[X_n^{2(k-1)}]^{1/2} \\
& \lesssim \mathbb{E}[(\Delta X_n)^2]^{1/2} n^{-(k-1)}
\end{aligned}
$$

and 

$$
\begin{aligned}
\mathbb{E}[(\Delta X_n)^2] &= \mathbb{E}[X_n^2] - 2 \mathbb{E}[X_nX_{n+1}] + \mathbb{E}[X_{n+1}^2] \\
&\leq 2\mathbb{E}[X_n^2] - 2 \mathbb{E}[X_n^2] + 2 \mathbb{E}[X_n\Delta X_n] \\
&\leq 2L^{-1} \kappa_d^{-1/d} \mathbb{E}[X_n^{2+1/d}] \\
&\leq 2L^{-1} \kappa_d^{-1/d} \mathbb{E}[X_n^3]^{\frac{2 + 1/d}{3}} \\
& \lesssim n^{-(2+1/d)}. 
\end{aligned}
$$

Consequently,

$$\frac{\mathbb{E}[X_n^k - X^k_{n+1}]}{\mathbb{E}[X_n^k]} \lesssim \frac{n^{-(2+1/d)/2 -(k-1)}}{n^{-k}} = n^{-1/(2d)} \rightarrow 0.$$

Let $\mu_k$ equal the limit of $n^k\mathbb{E}[X_n^k]$. To get convergence in distribution of $nX_n$, it suffices to prove that $$\limsup_{k\rightarrow \infty}\frac{\mu_{2k}^{1/2k}}{2k} < \infty$$ \cite[p.~109]{RD}. (\ref{mainlem:item:third}) gives $\mu_k \leq dm(B_{2L})k^k$ and therefore 

$$ \limsup_{k\rightarrow \infty}\frac{\mu_{2k}^{1/2k}}{2k} \leq \limsup_{k\rightarrow \infty}\frac{(dm(B_{2L})(2k)^{(2k)})^{1/2k}}{2k} = 1.$$

\end{proof}

\subsubsection{Proof of theorem \ref{wcpb}}\label{proofballs}

\begin{proof}
(\ref{wcpb:first}) The function \begin{equation*}G(x,\alpha)=(2L\omega_d)^{-1}\int_{|y|<|x|}|y|^{\alpha}|x-y|^{-(d-1)}dy\end{equation*} has the scaling property $G(\lambda x,\alpha)=\lambda^{\alpha+1} G(x,\alpha)$ for $\alpha \geq 0$. The scaling property implies \begin{equation}\label{rec}\mathbb{E}[X_n^\alpha]-\mathbb{E}[X_{n-1}^\alpha]=-c_{\alpha}\mathbb{E}[X_{n-1}^{\alpha+1}]\end{equation}where \begin{equation*}c_{\alpha}=(2L\omega_d)^{-1}\int_{|y|<1}(1-|y|^{\alpha})|u-y|^{-(d-1)}dy\end{equation*} and $u$ any unit vector. (\ref{exact2}) follows directly from the following recurrence relations (which follow from (\ref{rec})): \begin{equation}\label{rec_1}\sum_{k=0}^{n}\binom{n}{k}(-1)^{k}\mathbb{E}[X_{n-k}^j]=\triangle^n(X_0^j)=X_0^{j+n}(-1)^{n}\prod_{k=0}^{n}c_{k+j-1}.\end{equation}

(\ref{wcpb:monotone}) We use induction on $n$ for fixed $k$. Suppose $Z_n$ is like $X_n$ except that we start at $Z_0$ greater than $X_0$. The function $x \mapsto x^k-c_kx^{k+1}$ is increasing in $x$ for $x\leq L$ because $c_k \leq (2L)^{-1}$. Hence the statement is true for $n=1$. We have \begin{align*}\mathbb{E}[Z^k_{n+1}-X^k_{n+1}] &= \mathbb{E}[Z^k_n-X^k_n] -c_k\mathbb{E}[Z_n^{k+1}-X^{k+1}_n] \\ &\geq \mathbb{E}[Z^k_n-X^k_n]-c_k\mathbb{E}[(X_n+Z_n)(Z_n^k-X_n^k)] \\ &\geq \mathbb{E}[Z^k_n-X^k_n][1-c_k(X_0 + Z_0)] \\ &\geq 0  \end{align*} since $X_0 + Z_0 $ is less than or equal to $2L$. 

(\ref{wcpb:second}) For $d=1$, we have \begin{align*}n^j\mathbb{E}[X^j_n] &= n^j\sum_{k=0}^{n}\binom{n}{k}(-1)^kX_0^{k+j}(2L)^{-k}(j+k)^{-1} \\
&= (2L)^j n^j \int_{0}^{X_0/2L}t^{j-1}(1-t)^ndt \\
&\rightarrow (2L)^j j!\end{align*}as $n$ tends to infinity. This shows that $nX_n$ converges in distribution to an exponential with mean $2L$. 

(\ref{wcpb:third}) Suppose $\alpha \geq 1$. (\ref{rec}), the mean value theorem and Jensen's inequality yields
\begin{equation}
\mathbb{E}[X_n^{\alpha}]^{-1/\alpha} - X_0^{-1} \geq\alpha^{-1}c_{\alpha}\sum_{i=1}^n \mathbb{E}[X_{i-1}^{\alpha + 1}]\mathbb{E}[X_{i-1}^{\alpha}]^{-(1 + 1/\alpha)} \geq \alpha^{-1}c_{\alpha}n
\end{equation}which is equivalent to $\mathbb{E}[X_n^\alpha] \leq (X_0^{-1} + c_\alpha\alpha^{-1}n)^{-\alpha}$.

(\ref{wcpb:fourth}) For the lower bound, suppose $d>1$ and $R_n = X_n/X_{n-1}$ then $\mathbb{E}[X_n^{-1}] <\infty$ and 
\begin{align*}
\mathbb{E}[X_n^{-1}] - X_0^{-1} = \sum_{i=1}^n \mathbb{E}[X_{i-1}^{-1} \mathbb{E}[\frac{1-R_{i}}{R_{i}}|X_{i-1}]] = \sum_{i=1}^n \mathbb{E}[X_{i-1}^{-1} X_{i-1} \ell_d] = n\ell_d. 
\end{align*}Hence we obtain the desired lower bound $(X_0^{-1} + n\ell_d)^{-1} \leq \mathbb{E}[X_n^{-1}]^{-1} \leq \mathbb{E}[X_n]$.

(\ref{wcpb:sixth}) We first show that $n^k\mathbb{E}[X_n^k]$ converges for all $k\geq 1$. We have $$\frac{(n+1)^k\mathbb{E}[X_{n+1}^k]}{n^k\mathbb{E}[X_n^{k}]} \sim \frac{\mathbb{E}[X_{n+1}^k]}{\mathbb{E}[X_n^{k}]} = 1-c_{k}\frac{\mathbb{E}[X_n^{k+1}]}{\mathbb{E}[X_n^{k}]}
$$, $$\mathbb{E}[X_n^k] \geq \mathbb{E}[X_n]^k \gtrsim (n\ell_d)^{-k} $$ and $$\mathbb{E}[X_n^{k+1}] = O(n^{-(k+1)}).$$ Hence $$\lim_{n\rightarrow \infty} \frac{(n+1)^k\mathbb{E}[X_{n+1}^k]}{n^k\mathbb{E}[X_n^{k}]} = 1$$ and $n^k\mathbb{E}[X_n^{k}]$ converges. Set $\mu_k$ equal to the limit of $n^k\mathbb{E}[X_n^k]$. We have $$n^k\mathbb{E}[X_n^{k}] = n^kc_k\sum_{i=n}^{\infty}\mathbb{E}[X_i^{k+1}] \sim n^kc_k \mu_{k+1} \sum_{i=n}^{\infty}i^{-(k+1)} \sim c_k \mu_{k+1} k^{-1}.$$ Hence we have the recurrence relation $$\mu_k = c_k k^{-1} \mu_{k+1}$$ which implies that $$\mu_k = \mu_1 (k-1)! \prod_{i=1}^{k-1}c_i^{-1}.$$ To get convergence in distribution of $nX_n$, it suffices to prove that $$\limsup_{k\rightarrow \infty}\frac{\mu_{2k}^{1/2k}}{2k} < \infty$$ \cite[p.~109]{RD}. Since $c_i$ is increasing:$$\limsup_{k\rightarrow \infty}\frac{\mu_{2k}^{1/2k}}{2k} \leq \lim_{k \rightarrow \infty}\frac{(2k-1)c_{2k}^{-1}}{2k} = 2L\omega_d \left( \int_{|y| < 1} |u-y|^{-(d-1)} dy \right)^{-1}.$$ Hence $nX_n$ converges in distribution to a random variable $Y$ with moments $\mu_k$. We introduce the function $\varphi(x,L)$ which corresponds to $\mathbb{E}[Y]$ for starting parameters $X_0 = x$ and $L$. It follows directly from (\ref{wcpb:first}) that the function $\varphi$ has the following scaling property: $$\varphi(x,\lambda L) = \lambda \varphi(\lambda^{-1}x,L).$$ Now if the limiting distribution doesn't depend on the starting value $X_0=x$ then $\varphi(\lambda^{-1}x,L) = \varphi(x,L)$ and we have the scaling property $$\varphi(x,\lambda L) = \lambda \varphi(x,L).$$ Finally, to find $\mathbb{E}[Y]$ for a fixed $L$, we use $$\mathbb{E}[Y] = \varphi(x,L) = \varphi(L,L) = L\varphi(1,1).$$ But we know that $$\varphi(1,1) = \lim_{n\rightarrow \infty} n \sum_{k=0}^{n}\binom{n}{k}(-2)^k \prod_{i=0}^{k}\overline{c}_i $$ where $$\overline{c}_k = \omega_d^{-1}\int_{|y|<1}(1-|y|^k)|u-y|^{-(d-1)} \, dy$$ for $k>0$ and $\overline{c}_0 $ is set to $1$. To complete the proof, we need to prove that the limiting distribution doesn't depend on $X_0$. Define $X_n$ and $Z_n$ as in the proof of (\ref{wcpb:monotone}). It follows from the monotonicity property (\ref{wcpb:monotone}) that \begin{align*}\mathbb{E}[Z_n] &\leq  \mathbb{P}(Z_N \leq X_0) \mathbb{E}[X_{n-N}] +  \mathbb{P}(Z_N > X_0) X_0 \\ & \leq  \mathbb{E}[X_{n-N}] + \mathbb{P}(Z_N > X_0) X_0 \\ &\leq \mathbb{E}[X_{n-N}] + \mathbb{E}[Z_N^2]X_0^{-1}\end{align*} for all $n > N$. Let $N = \lfloor n^{2/3} \rfloor$. Given that $\mathbb{E}[Z_N^2] = O(N^{-2})$, and letting $L_X = \lim_{n \rightarrow \infty} n\mathbb{E}[X_n]$ and $L_Z = \lim_{n \rightarrow \infty} n\mathbb{E}[Z_n]$, we have:

\begin{align*}
\lim_{n \rightarrow \infty} n\mathbb{E}[Z_n] 
&\leq \lim_{n \rightarrow \infty} \left( \frac{n}{n-N} \right) (n-N)\mathbb{E}[X_{n-N}] + \lim_{n \rightarrow \infty} \left( \frac{n}{N} \right) O(N^{-2})X_0^{-1} \\
&= \left( \lim_{n \rightarrow \infty} \frac{1}{1-n^{-1/3}} \right) \cdot L_X + \lim_{n \rightarrow \infty} O(n \cdot n^{-4/3})X_0^{-1} \\
&= (1) \cdot L_X + \lim_{n \rightarrow \infty} O(n^{-1/3})X_0^{-1} \\
&= L_X \\
&= \lim_{n \rightarrow \infty} n\mathbb{E}[X_n]
\end{align*}

(\ref{wcpb:fifth}) The volume of the symmetric difference between two balls of radius $r$ with centers a distance $x$ apart equals  
\begin{equation}\label{ballsymmd}
4\kappa_{d-1}r^{d}\int_{\theta(x)}^{\pi/2} \sin^d(t) dt,\; \theta(x) = \arccos(x/2r \wedge 1).
\end{equation}We introduce the function $\phi(x) = \int_{\theta(x)}^{\pi/2} \sin^d(t) dt$ with $x$ between $-1$ and $1$. On the interval $[-\cos(\theta_0),\cos(\theta_0)]$ with $0 < \theta_0 < \pi/2$ , the second derivative is bounded and $\phi'(0) = 1$; consequently, we can write $\phi(x) = x + o(|x|)$ for $|x| \leq \cos(\theta_0)$. Applying (\ref{ballsymmd}) to the random sequence $X_n$:
\begin{align*}
nm(A_n \triangle A^*) &= 4\kappa_{d-1}r^{d}n\phi(X_n/2r \wedge 1) \\
&=  4\kappa_{d-1}r^{d}n(X_n/2r \wedge \cos(\theta_0)) + no(X_n/2r \wedge \cos(\theta_0)) \\
&\Rightarrow 4\kappa_{d-1}r^{d} Y
\end{align*}.

where $\Rightarrow$ refers to convergence in distribution. 

\end{proof}

\textbf{Remark:} For $d=3$, the constants $\overline{c}_k$ can be computed exactly: $$\overline{c}_k = \begin{cases}\frac{1}{2}(1-\frac{2\log(2) + 2\sum_{i=1}^{k/2}\frac{1}{2i}}{k+1}), \, k \, \textrm{even}\\ \frac{1}{2}(1-\frac{2\sum_{i=1}^{(k+1)/2}\frac{1}{2i-1}}{k+1}), \, k \,\textrm{odd} \end{cases}$$

\section{Random polarizations as a non-homogeneous Markov chain}

We introduce another sequence of random polarizations where the transition probability depends on the current state of the Markov chain. Consider a set $A$ of finite measure in $\mathbb{R}^d$. Define the following random sequences:

\begin{itemize}

\item $A_n = A^{\sigma_1\cdots\sigma_n}$
\item $X_n = m(A_n/A^*)$
\item $\pi_n(x) = X_n - m\left(A_n^{\sigma_{0,x}}/A^*\right)$
\item $g_n(x) \in L^1(\mathbb{R}^d)$ .
\end{itemize} 

\noindent The probability distribution for $\sigma_n$ is as follows:

\begin{equation}
\mathbb{P}(\sigma_{n+1}(0) \in A \, | \, X_n) = \frac{\int_{A} g_n(x)\pi_n(x) \, dx}{\int_{\mathbb{R}^d}g_n(x) \pi_n(x) \, dx}.
\end{equation}

\begin{Lem}
\begin{equation}
\int_{\mathbb{R}^d} g_n(t) \pi_n(t) \, dt = \frac{1}{\sqrt{2}} \int_{A_n/A^*} \int_{A^*/A_n} g_n(u(x,y)) |J_d u(x,y)|\, dx \, dy
\end{equation}

where $J_du$ is the $d$-dimensional Jacobian of $u$ with determinant $|J_d u|$.

\end{Lem}

\begin{proof}

By the coarea formula: 

\begin{equation} \int_{A_n/A^*} \int_{A^*/A_n} g(u(x,y))|J_d u(x,y)|\, dx \, dy = \int_{\mathbb{R}^d} g(t)H_{d}(\{u=t\})\,dt \end{equation}

where $H_d$ is the $d$-dimensional Hausdorff measure. The preimage of $t$ under $u$ is the graph 

$$\{(x,\sigma_{0,t}(x))\ \,|\, x \in E_{n,t} \}$$

where $E_{n,t} = A_n/A^* \cap \sigma_{0,t}^{-1}(A^*/A_n)$. The $d$-dimensional Hausdorff measure of this graph equals 

$$\int_{E_{n,t}}\sqrt{1 + |\nabla \sigma |^2} dx = \sqrt{2} m(E_{n,t}) = \sqrt{2} \pi_n(t)$$ where we have used the fact that $|\nabla \sigma| = 1$ since $\sigma$ is an isometry.

\end{proof}

\begin{Lem}
\begin{equation}
|J_d u(x,y)| = \left(\frac{\left||x|^2-|y|^2 \right|}{|x-y|^2}\right)^{d-1} \sqrt{2}
\end{equation}

\end{Lem}

\begin{proof}

We have \begin{equation}\label{fast-jac-limit} \lim_{\delta \rightarrow 0 } \delta^{-2d} \kappa_d^{-2} \int_{B_{x_0,\delta}} \int_{B_{y_0,\delta}}  |J_d u(x,y)| \, dx \, dy = |J_d u(x_0,y_0)|.\end{equation} Following the proof of the previous lemma, the integral on the left-side also equals

\begin{equation}\label{fast-jac-int} \sqrt{2} \int_{u(B_{x_0,\delta} \times B_{y_0,\delta})}  m(E_{t}) \, dt \end{equation}

where $$E_t = B_{y_0,\delta} \cap \sigma_{0,t}( B_{x_0,\delta}) = B_{y_0,\delta} \cap B_{\sigma_{0,t}(x_0),\delta}. $$ We have $$m(E_t) = 4 \delta^d \kappa_{d-1} \phi\left(\frac{|\sigma_{0,t}(x_0) - y_0|}{2\delta}\right)$$ where

$$\phi(x) = \int_0^{\arccos(x)}\sin^d(t) \, dt.$$

Making the change of variable $y = (\sigma_{0,t}(x_0)-y_0)/(2\delta)$, the integral in \ref{fast-jac-int} equals

\begin{equation*} 
2^{d+2}\delta^{2d} \kappa_{d-1} \int_{B_1} \left(\frac{\left||x_0|^2-|2\delta y+y_0|^2 \right|}{|x_0-y_0 - 2\delta y|^2}\right)^{d-1} \phi(|y|) \, dy 
\end{equation*} and the limit \ref{fast-jac-limit} equals 

\begin{equation}\frac{ 2^{d + 2 -\frac{1}{2}}\kappa_{d-1}}{\kappa_d^2}\left(\frac{\left||x_0|^2-|y_0|^2 \right|}{|x_0 - y_0|^2}\right)^{d-1} \int_{B_1} \phi(|y|) \, dy.  \end{equation} 

By Fubini's theorem: 

\begin{equation}
\int_{B_1} \phi(|y|) \, dy = \kappa_d \int_0^{\frac{\pi}{2}} \sin^d(t) \cos^d(t) \, dt = \frac{\kappa_d \Gamma(\frac{d}{2} + \frac{1}{2})^2}{2d!}.
\end{equation}

Using the identity $$ \frac{\kappa_d}{\kappa_{d-1}} = \frac{\sqrt{\pi} \Gamma(\frac{d}{2} + \frac{1}{2})}{\Gamma(\frac{d}{2} + 1)}$$  and the duplication identity

$$\Gamma \left(\frac{d}{2}+ 1 \right) \Gamma\left(\frac{d}{2}+ \frac{1}{2}\right) = 2^{-d} \sqrt{\pi} d!$$

gives the desired result.

\end{proof}

\begin{Thm}

\label{fast-mainthm-lbound}

Suppose $A \subset B_L$ and set $g_n(x) = |x|^{-(d-1)}$ for all $n$. If $E_n = u(A_n/A^* \times A^*/A_n)$ and $$Y_n = \left(\frac{m(E_n)}{k_d}\right)^{1/d}$$ then 
\begin{equation}
\mathbb{E}[X_n - X_{n+1} | X_n] \geq \frac{\sqrt{2} X_n^2}{ \operatorname{Per}(B_{2L}) \cdot Y_n}.
\end{equation}

\end{Thm}

\begin{proof}

By Jensen's inequality:

$$\mathbb{E}[X_{n} - X_{n+1} | X_n] = \frac{\int_{E_n} g_n(x)\pi_n^2(x) \, dx}{\int_{E_n} g_n(x) \pi_n(x) \, dx} \geq  \frac{\int_{E_n} g_n(x)\pi_n(x) \, dx}{\int_{E_n} g_n(x) \, dx}.$$

By the previous lemmas, the integral $\int_{E_n} g_n(x)\pi_n(x) \, dx$ equals

$$\sqrt{2}\int_{A^*/A_n} \int_{A_n/A^*} |x-y|^{-(d-1)}\, dx \, dy. $$ 

The integral above is bounded below by $\sqrt{2}(2L)^{-(d-1)}X_n^2$. The Riesz rearrangement inequality gives

$$\int_{E_n} g_n(x) \, dx \leq \omega_d Y_n$$ where $$Y_n = \left(\frac{m(E_n)}{k_d}\right)^{1/d}, \, E_n = u(A_n/A^* \times A^*/A_n).$$ 

\end{proof}

\begin{Lem}
\label{lem:exponential-tail-decay}
Let $A_n$ be a sequence of measurable sets generated by successive random polarizations of a measurable set $A\subset B_L \subset \mathbb{R}^d$ according to the probability measure $\mathbb{P}$ defined in \eqref{eq:standard-prob}. 

For a fixed threshold radius $t > R^*$ where $R^*$ is the radius of $A^*$, we define the tail set $E_n = A_n/A^* \cap \{|x| > t\}$. Then, the expected measure of the tail set decays exponentially:
\[
\mathbb{E}[m(E_n)] \le m(E_0) \exp\left( -\rho_d \cdot n \right),
\]
where the explicit structural contractive constant $\rho_d$ is given by:
\[
\rho_d = \left( \frac{t}{2L} \right)^d - \left( \frac{R^*}{2L} \right)^d > 0.
\]

Similarly, for a fixed threshold radius $0 < t < R^*$, we define the inner tail set $E_n = (A^* \setminus A_n) \cap \{|x| < t\}$. Then, the expected measure of this inner tail set decays exponentially:
\begin{equation}
\mathbb{E}[m(E_n)] \le m(E_0) \exp\left( -\rho_d \cdot n \right),
\end{equation}
where the explicit structural contractive constant $\rho_d$ is given by:
\begin{equation}
\rho_d = \left( \frac{R^*}{2L} \right)^d - \left( \frac{t}{2L} \right)^d > 0.
\end{equation}

\end{Lem}

\begin{proof}
We prove the first statement only as the proof for the second statement is exactly the same. The expected volume step drop under this polarization mechanism satisfies:

\begin{align}
\mathbb{E}[m(E_n) - m(E_{n+1}) \mid A_n] &\ge (2L \omega_d)^{-1} (2L)^{1-d} m(F_n) m(E_n) \nonumber \\
&= (2L\omega_d)^{-1} (2L)^{1-d} \Big[ \omega_d(t^d - (R^*)^d) + m(E_n) \Big] m(E_n), \label{eq:decoupled-step}
\end{align}
where we have utilized the identity $m(F_n) = \omega_d(t^d - (R^*)^d) + m(E_n)$. Since $m(E_n) \ge 0$, dropping the tail term yields a conservative, static lower bound on the sink capacity. Simplifying the coefficients then directly yields:
\begin{align*}
\mathbb{E}[m(E_n) - m(E_{n+1}) \mid A_n] &\ge \frac{1}{2L \omega_d} (2L)^{1-d} \omega_d \left( t^d - (R^*)^d \right) m(E_n) \\
&= \left[ \left( \frac{t}{2L} \right)^d - \left( \frac{R^*}{2L} \right)^d \right] m(E_n).
\end{align*}
Setting $\rho_d = \left( \frac{t}{2L} \right)^d - \left( \frac{R^*}{2L} \right)^d$, this relation simplifies to the linear contraction:
\[
\mathbb{E}[m(E_{n+1}) \mid A_n] \le (1 - \rho_d) m(E_n).
\]
Taking the total expectation across both sides by applying the tower property of conditional expectation, we iterate the discrete system inductively from $n=0$:
\[
\mathbb{E}[m(E_n)] \le m(E_0) (1 - \rho_d)^n.
\]
Applying the standard inequality $1 - \rho_d \le \exp(-\rho_d)$ completes the proof.
\end{proof}

For $d>1$, we define the sets $F_n$ and $G_n$ as follows:
$$
\begin{aligned}
F_n &= \left\{ (x,y) \in (A_n \setminus A^*) \times (A^* \setminus A_n) : |x| \le r + \delta_n, \, |y| > r - \delta_n \right\} \cap u^{-1}(B_{\lambda_n}^c) \\
G_n &= \left\{ (x,y) \in (A_n \setminus A^*) \times (A^* \setminus A_n) : |x| \le r + \delta_n, \, |y| > r - \delta_n \right\} \cap u^{-1}(B_{\lambda_n})
\end{aligned}
$$
The function $g_n(x)$ is given by:
$$
g_n(x) = \begin{cases}
|x|^{-(d-1)} \mathbbm{1}_{B_{\lambda_n}^c}, & \text{for } m(F_n) > m(G_n) \text{ and } m(F_n \cup G_n) \ge \frac{1}{2}X_n^2 \\
|x|^{-(d-1)} \mathbbm{1}_{B_{\lambda_n}}, & \text{for } m(F_n) \le m(G_n) \text{ and } m(F_n \cup G_n) \ge \frac{1}{2}X_n^2 \\
|x|^{-(d-1)} \pi_n(x)^{-1}, & \text{otherwise}
\end{cases}
$$

The sequence $\delta_n$ and $\lambda_n$ will be defined explicitly in the proof of the following theorem.

\begin{Thm}
\label{fast-maincor-2}

For $d>1$, $\mathbb{E}[X_n] = \mathcal{O}\left( n^{-\left(2 - \frac{1}{d}\right)} (\log n)^{1 - \frac{1}{d}} \right)$.

\end{Thm}

\begin{proof}

If $m(F_i \cup G_i) \le \frac{1}{2}X_i^2$ is true more than $\left\lfloor \frac{n}{2} \right\rfloor$ times after $n$ random polarizations governed by the sequence $g_n$, then by Lemma \ref{lem:exponential-tail-decay}: 

\begin{equation}\mathbb{E}[X_n] \leq \frac{\sqrt{2}}{\sqrt{2}-1}  X_0 \exp\left( -\rho_d \left\lfloor \frac{n}{2} \right\rfloor \right) \end{equation}

where \begin{equation}
\rho_d = \left |\left( \frac{1 \pm \delta_n}{2L} \right)^d - \left( \frac{1}{2L} \right)^d \right |. 
\end{equation}

Now we consider the case when $m(F_i \cup G_i) < \frac{1}{2}X_i^2$ more than $\left\lfloor \frac{n}{2} \right\rfloor$ times

We first consider the case $m(F_n \cup G_n) \ge \frac{1}{2}X_n^2$. If $m(F_n) > m(G_n)$ then

\begin{align}
\mathbb{E}[X_{n} - X_{n+1} | X_n]  &\geq  \frac{\sqrt{2}\int\int \mathbbm{1}_{F_n} |x-y|^{-(d-1)} \, dx \, dy}{(2L)\omega_d } \\
& \geq \frac{(2L)^{-(d-1)}\lambda_n^{d-1} m(F_n)}{2^{d - \frac{1}{2}} \omega_d L \delta_n^{d-1}} \\
& \geq \frac{\lambda_n^{d-1} X_n^2}{2^{2d + \frac{1}{2}} \omega_d L^d \delta_n^{d-1}}.
\end{align}

And similarly if $m(G_n) \geq m(F_n)$ then 

\begin{align}
\mathbb{E}[X_{n} - X_{n+1} | X_n]  &\geq  \frac{\sqrt{2} \int\int \mathbbm{1}_{G_n} |x-y|^{-(d-1)} \, dx \, dy}{\lambda_n} \\
&\geq \frac{X_n^2}{2^{d - \frac{1}{2}}L^{d-1}\lambda_n}. 
\end{align}

In particular, if we set 

\begin{equation}
\lambda_n = \left( 2^{d+1} \omega_d L \right)^{\frac{1}{d}} \delta_n^{1 - \frac{1}{d}}
\end{equation}

then 

\begin{equation}
\mathbb{E}[X_{n} - X_{n+1} \mid X_n] \geq \frac{X_n^2}{2^{d-\frac{1}{2}}L^{d-1}\lambda_n} = \frac{X_n^2}{C_d \delta_n^{1 - \frac{1}{d}}}
\end{equation}
where $C_d = 2^{d - \frac{1}{2}} L^{d-1} \left( 2^{d+1} \omega_d L \right)^{\frac{1}{d}}$. We then obtain the upper bound:

\begin{equation}
\mathbb{E}[X_n] \leq \left( X_0^{-1} + \sum_{k=\lfloor n/2 \rfloor}^{n-1} \frac{1}{C_d \delta_k^{1 - \frac{1}{d}}} \right)^{-1}
\end{equation}

So combining both scenarios, we finally obtain the following upper bound: 

\begin{equation}
\mathbb{E}[X_n] \leq \max \left\{ 2 \sqrt{X_0} \exp\left( -\frac{\rho_d}{2} \left\lfloor \frac{n}{2} \right\rfloor \right), \, \left( X_0^{-1} + \sum_{k=\lfloor n/2 \rfloor}^{n-1} \frac{1}{C_d \delta_k^{1 - \frac{1}{d}}} \right)^{-1} \right\}
\end{equation}

To complete the proof, we evaluate the asymptotic behavior of each term inside the maximum under the choice of $\delta_n = c_d \frac{\log n}{n}$, where $c_d = \frac{8(2L)^d(2d-1)}{d^2}$.

We begin by estimating the contractive constant $\rho_d$. For large $n$, we have $\delta_n \to 0$. Applying a first-order Taylor expansion yields:
\begin{equation}
\rho_d = \frac{d}{(2L)^d} \delta_n + \mathcal{O}(\delta_n^2).
\end{equation}
Substituting our explicit choice of $\delta_n$ into this expansion gives:
\begin{equation}
\rho_d \sim \frac{d}{(2L)^d} \cdot \left( \frac{8(2L)^d (2d - 1)}{d^2} \right) \frac{\log n}{n} = \frac{8(2d-1)}{d} \frac{\log n}{n}.
\end{equation}
Plugging this expression into the first branch of the maximum bound, we find:
\begin{align}
 \exp\left( -\frac{\rho_d}{2} \left\lfloor \frac{n}{2} \right\rfloor \right) &\le \exp\left( -\frac{4(2d-1)}{d} \frac{\log n}{n} \cdot \left(\frac{n}{2} - 1\right) \right) \nonumber \\
&\sim  \exp\left( -\frac{2(2d-1)}{d} \log n \right) \nonumber \\
&=  n^{-\left(4 - \frac{2}{d}\right)}.
\end{align}

Now we estimate the second term. Substituting $\delta_k$ into the sum over the tail window $k = \lfloor n/2 \rfloor, \dots, n-1$, the general term becomes:
\begin{equation}
\frac{1}{C_d \delta_k^{1 - \frac{1}{d}}} = \frac{1}{C_d c_d^{1 - \frac{1}{d}}} \left( \frac{k}{\log k} \right)^{1 - \frac{1}{d}}.
\end{equation}
Approximating the sum by an integral over the index range yields:
\begin{align}
\sum_{k=\lfloor n/2 \rfloor}^{n-1} \left( \frac{k}{\log k} \right)^{1 - \frac{1}{d}} &\sim \int_{n/2}^{n} \left( \frac{x}{\log x} \right)^{1 - \frac{1}{d}} \, dx \nonumber \\
&\sim \frac{1}{(\log n)^{1 - \frac{1}{d}}} \int_{n/2}^{n} x^{1 - \frac{1}{d}} \, dx \nonumber \\
&= \frac{1}{(\log n)^{1 - \frac{1}{d}}} \cdot \frac{1}{2 - \frac{1}{d}} \left( n^{2 - \frac{1}{d}} - \left(\frac{n}{2}\right)^{2 - \frac{1}{d}} \right) \nonumber \\
&= \frac{1 - 2^{-\left(2 - \frac{1}{d}\right)}}{2 - \frac{1}{d}} \cdot \frac{n^{2 - \frac{1}{d}}}{(\log n)^{1 - \frac{1}{d}}}.
\end{align}
Inverting this cumulative sum dominates the initial value $X_0^{-1}$ as $n \to \infty$, which gives the asymptotic behavior of the second branch:
\begin{equation}
\left( X_0^{-1} + \sum_{k=\lfloor n/2 \rfloor}^{n-1} \frac{1}{C_d \delta_k^{1 - \frac{1}{d}}} \right)^{-1} \sim \left( \frac{C_d c_d^{1 - \frac{1}{d}} \left(2 - \frac{1}{d}\right)}{1 - 2^{-\left(2 - \frac{1}{d}\right)}} \right) n^{-\left(2 - \frac{1}{d}\right)} (\log n)^{1 - \frac{1}{d}}.
\end{equation}

Comparing the polynomial exponents of both branches, we observe that:
\begin{equation}
4 - \frac{2}{d} = 2\left(2 - \frac{1}{d}\right) > 2 - \frac{1}{d}.
\end{equation}
Consequently, the exponential tail decay term decays strictly faster than the harmonic sum term. The maximum is therefore asymptotically dominated by the second branch, establishing the final baseline convergence rate:
\begin{equation}
\mathbb{E}[X_n] = \mathcal{O}\left( n^{-\left(2 - \frac{1}{d}\right)} (\log n)^{1 - \frac{1}{d}} \right).
\end{equation}

\end{proof}

The following proposition shows that there exists a class of measurable sets in $\mathbb{R}$ that admit polarizations converging to their corresponding Schwarz symmetrization at an exponential rate.

\begin{Lem}If $d=1$ and $A = \{f>t\}$ for a non-negative Lipschitz continuous function $f$, with Lipschitz constant $C$ and distribution function that is differentiable at $t$, then there exists a sequence of polarizations $A_n$ of $A$ such that \begin{displaymath}m(A_n \triangle A^*) \leq (1-C^{-1}(-\mu_f'(t)+C^{-1})^{-1})^n\end{displaymath} \end{Lem}

for $n$ sufficiently large. 

\begin{proof}

Suppose $f$ is Lipschitz continuous with Lipschitz constant $C$. If $f(x_0) > t + \lambda $ and $f(y_0) \leq t - \lambda$, then, by the triangle inequality, $f(x) > t$ and $f(y) \leq t$ whenever $|x-x_0| \leq C^{-1}\lambda$ and $|y-y_0| \leq C^{-1}\lambda$. If we add the conditions $|x_0| > r_f(t) + C^{-1}\lambda$ and $|y_0| \leq r_f(t) - C^{-1}\lambda$, then $|x| >r_f(t) $, $|y| \leq r_f(t)$, $x \in A/A^*$, $y \in A^*/A$, and $$m(A^{\sigma_{x_0,y_0}} \triangle A^*) \leq m(A \triangle A^*) -4 C^{-1}\lambda.$$ Note that $$A/A^* \cap \{f \leq t+\lambda\} = \{x: t < f(x) \leq t + \lambda, \, |x| >r_f(t) \}$$ and $$A^*/A \cap \{f > t -\lambda \} = \{y: t-\lambda < f(y) \leq t, \, |y| \leq r_f(t) \}.$$ Since the distribution function of $f$ is differentiable at $t$:

$$m(A/A^* \cap \{f \leq t+\lambda\}) \leq \mu_f(t)-\mu_f(t+\lambda)  = -\mu_f'(t)\lambda (1 + o(1)) $$

and $$m(A^*/A \cap \{f > t -\lambda \}) \leq \mu_f(t-\lambda) - \mu_f(t) = -\mu_f'(t)\lambda (1+o(1))$$

as $\lambda \rightarrow 0$. Similarly, we have 

$$m(A/A^* \cap B_{r_f(t) + C^{-1}\lambda}) \leq m(B_{r_f(t) + C^{-1}\lambda}) - m(B_{r_f(t)}) = 2C^{-1}\lambda $$

and $$ m((A^*/A) / B_{r_f(t) - C^{-1}\lambda}) \leq m(B_{r_f(t)}) - m(B_{r_f(t) - C^{-1}\lambda}) = 2C^{-1}\lambda. $$

Setting $$\lambda = \frac{1}{4} (-\mu_f'(t)+C^{-1})^{-1}m(A \triangle A^*)$$, we see from the inequalities above that there exists points $x_0$ and $y_0$ such that 

$$m(A^{\sigma_{x_0,y_0}} \triangle A^*) \leq m(A \triangle A^*)(1-C^{-1}(-\mu_f'(t)+C^{-1})^{-1}) $$

provided that $m(A \triangle A^*)$ is less than $\delta_f$ where $\delta_f$ depends only on $\mu_f$. The proof of the proposition is completed by noting that if $f$ is Lipschitz continuous with Lipschitz constant $C$ then any polarization of $f$ is also Lipschitz continuous with Lipschitz constant less than or equal to $C$. 

\end{proof}

\begin{Thm}

For every compact set $A$ in $B_L$ with finite Minkowski perimeter there exists a sequence of polarizations $A_n$ of $A$ such that 

\begin{displaymath}m(A_n \triangle A^*) \leq \left(\frac{\operatorname{Per}(A)}{1 + \operatorname{Per}(A)}\right)^n\end{displaymath}

for $n$ sufficiently large. 

\end{Thm}

\begin{proof}

We can express the interior of $A$ as the disjoint union of open intervals: $$A^{\circ} = \bigcup_{k=-\infty}^{\infty}(a_k,b_k). $$ Then we define the Lipschitz function $f$: 

\begin{equation}
f(x) = \frac{L}{2} + \frac{a_k + b_k}{2} - |x - \frac{a_k + b_k}{2}|, \, x \in [\frac{b_{k-1} + a_k}{2},\frac{b_k + a_{k+1}}{2}]. 
\end{equation}

The Lipschitz constant is 1 and $\{f > L/2\} = A^{\circ}$. Furthermore, we have 

$-\mu^{'}_f(L/2) = \operatorname{Per}(A).$ Applying the previous lemma completes the proof. 

\end{proof}

\bibliography{roc_polarization.bbl}

\providecommand{\bysame}{\leavevmode\hbox to3em{\hrulefill}\thinspace}
\providecommand{\MR}{\relax\ifhmode\unskip\space\fi MR }
\providecommand{\MRhref}[2]{%
  \href{http://www.ams.org/mathscinet-getitem?mr=#1}{#2}
}
\providecommand{\href}[2]{#2}
\begin{thebibliography}{1}

\bibitem{AB}
Almut Burchard, \emph{Rate of convergence of random polarizations}, Preprint
  arXiv:1108.5500 (2011), 5 Pages, \url{http://arxiv.org/abs/1108.5500}.

\bibitem{ABQD}
Almut Burchard and Qin Deng, \emph{On the rate of convergence of random
  polarizations on the sphere},  (2020), 10 Pages,
  \url{https://www.math.toronto.edu/dengqin/polarization.pdf}.

\bibitem{AMF}
Almut Burchard and Marc Fortier, \emph{Random polarizations}, Preprint
  arXiv:1104.4103v4 (2012), 26 Pages, \url{http://arxiv.org/abs/1104.4103}.

\bibitem{RD}
Rick Durrett, \emph{Probability: Theory and examples}, Duxbury Advanced Series,
  2005.

\end{thebibliography}

\end{document}